\documentclass[a4paper]{amsart}

\usepackage{amssymb}
\usepackage{amsthm}
\usepackage{amsmath}
\usepackage{mathrsfs}
\usepackage[colorlinks, allcolors=blue]{hyperref}

\newtheorem{theorem}{Theorem}
\newtheorem{theoremA}{Theorem}

\newtheorem{lemma}{Lemma} 

\theoremstyle{definition} 

\newtheorem*{problem*}{Problem}
\newtheorem{remark}{Remark}

\newcommand{\ind}{\operatorname{ind}}

\newcommand{\bR}{\mathbb R}

\begin{document}

{
\renewcommand*{\thefootnote}{$\star$}

\title[There are No Product and Subgroup Theorems]{There are No Product and Subgroup Theorems\\
for the Covering Dimension\\ of Topological Groups}
\footnotetext[0]{This work was financially supported by the Russian Science Foundation, grant~22-11-00075-P.} 

}

\author{Ol'ga Sipacheva}

\address{Faculty of Mechanics and  Mathematics, 
M.~V.~Lomonosov Moscow State University, Moscow, Russia}

\email{ovsipa@gmail.com}

\begin{abstract}
Strongly zero-dimensional topological groups $G_1$, $G_2$, and $G$ such that $G_1\times G_2$ has 
positive covering dimension and $G$ contains a closed subgroup of positive covering dimension are 
constructed. Moreover, all finite powers of $G_1$ are Lindel\"of and $G_2$ is 
second-countable. An example of a strongly zero-dimensional space $X$ whose free, free Abelian, and free Boolean 
topological groups have positive covering dimension is also given. 
\end{abstract}



\maketitle

This paper is concerned with the covering dimension of topological groups. There are two definitions of covering 
dimension, in the sense of \v Cech and in the sense of Kat\v etov; to differentiate them, following \cite{Ch}, we 
denote the former by $\dim$ and the latter by $\dim_0$. 

In 1989 Shakhmatov asked whether the inequality $\dim_0(G \times H) \le \dim_0 G + \dim_0 H$ holds for arbitrary 
topological groups $G$ and $H$ and proved that the answer is positive for precompact groups\footnote{In 
\cite{88}, as well as in most papers cited below, Kat\v etov covering dimension is denoted by $\dim$ rather than 
by $\dim_0$.} \cite{88}. Various versions of this question can be found in \cite{A-vM}. We construct two 
topological groups $G$ and $H$ such that all finite powers of $G$ are Lindel\"of, $H$ is second-countable, and 
$\dim_0(G\times H)>\dim_0 G + \dim_0 H =0$, thereby answering (in the negative) Shakhmatov's question and 
Questions~6.9 and~6.14 of \cite{A-vM}. A modification of this example gives a negative answer to Arkhangel'skii's 
old question of whether the free (free Abelian) topological group of any strongly zero-dimensional space is 
strongly zero-dimensional (see \cite[p.~964]{Arh89} and \cite[Problem~8.17]{A-vM}).

In the same paper \cite{88} Shakhmatov also asked whether the inequality $\dim_0 H\le \dim_0 G$ holds for an 
arbitrary subgroup $H$ of an arbitrary topological group $G$ (see also \cite[Problem~6.9]{Tk4} and 
\cite[Question~6.1]{A-vM}). He proved that the answer is positive if $G$ is a locally pseudocompact or 
Lindel\"of $\Sigma$ group~\cite{88, Sh}. It is also known that if $H$ is $\bR$-factorizable, then $\dim_0 H\le 
\dim_0 G$~\cite[Theorem~2.7]{Tk}. In particular, if $H$ is Lindel\"of, then $\dim H=\dim_0 H\le \dim_0 G$. In 
this paper we construct a topological group $G$ with $\dim_0 G=0$ which contains a closed subgroup $H$ of 
positive covering dimensions $\dim_0$ and $\dim$. Moreover, $H$ is the product of two Lindel\"of groups, one of 
which is second-countable.

\section{Preliminaries}

For convenience, we assume all topological spaces and groups considered in this paper to be Tychonoff. 

Suppose given a set $X$ and a family $\mathscr F$ of its subsets. If there exists an integer $n\ge 
-1$ such that every point of $X$ belongs to at most $n+1$ elements of $\mathscr F$, then the smallest such 
$n$ is called the \emph{order} of $\mathscr F$; otherwise, the order of $\mathscr F$ is infinite. Given a 
topological space $X$, the \emph{\v Cech covering dimension $\dim X$} of $X$ is the smallest integer $n$ for 
which any finite open cover of $X$ has a finite open refinement of order $n$, provided that such an integer 
exists; if it does not exist, then $\dim X$ is infinite. The definition of the \emph{Kat\v etov covering 
dimension $\dim_0 X$} is similar but uses covers by cozero sets: $\dim_0 X$ equals the smallest integer $n$ for 
which any finite cozero cover of $X$ has a finite cozero refinement of order $n$ if such an integer exists and is 
infinite otherwise. For normal spaces, the dimensions $\dim$ and $\dim_0$ coincide (see, e.g., 
\cite[Proposition~11.2]{Ch}). A space $X$ for which $\dim_0 X =0$ is said to be \emph{strongly zero-dimensional}. 

In what follows, we mention the \emph{small inductive dimension} $\ind X$ of a space $X$; we refer the reader 
to \cite{Ch} or \cite{Eng} for its definition, because it is not important for our purposes. It is only 
important that $\ind X=0$ if and only if $X$ has a base consisting of clopen sets. A space $X$ with $\ind X=0$ is 
said to be \emph{zero-dimensional}. Obviously, any strongly zero-dimensional space is zero-dimensional. 

A state-of-the-art presentation of the dimension theory of topological spaces is given in the 
highly recommended book~\cite{Ch} of Michael Charalambous. 

\begin{remark}
\label{r1}
 There are many examples of spaces $X$ for which $\dim X> \dim_0 X$ (see, e.g., \cite{Ch}), but the author is 
unaware of any example of a space $X$ for which $\dim_0 X> \dim X$. However, even if such examples exist, we 
have $\dim X = 0$ whenever $\dim_0 X =0$, because any finite disjoint open cover of $X$ consists of clopen 
sets, which are obviously cozero. 
\end{remark}

We use the notation $\mathbb R$ for the set of real numbers, $\mathbb N$ for the set of positive integers, 
and $\omega$ for the set of nonnegative integers. By $\oplus$ we denote the topological sum of spaces and     
by $|A|$, the cardinality of a set $A$. 

A subset $Y$ of a space $X$ is said to be \emph{$C$-embedded} in $X$ if any real-valued 
continuous function on $Y$ has a continuous extension to $X$, and $Y$ is \emph{$z$-embedded} in $X$ if every zero 
set of $Y$ is the trace on $Y$ of some zero set of $X$. A topological space admitting a coarser metrizable 
topology is said to be \emph{submetrizable}. 

To distinguish groups without topology from topological groups, we refer to the former as \emph{abstract 
groups}. 

The topology of a topological group $G$ is \emph{linear} if the open subgroups of $G$ form a base of 
neighborhoods of the identity element in~$G$. Clearly, all groups with linear topology are zero-dimensional, 
because any open subgroup of any topological group is closed. 

A \emph{Boolean group} is a group in which all elements are of order 2. All such groups are Abelian; moreover, 
all of them are vector spaces of the two-element field $\mathbb F_2$ and hence are free. Given a set $X$, the 
Boolean group $B(X)$ with basis $X$ is nothing but the set $[X]^{<\omega}$ of finite subsets of $X$ endowed with 
the operation of symmetric difference, which  plays the role of addition. The zero element is the empty set. 
Each point $x\in X$ is identified with the singleton $\{x\}$. 

For a Tychonoff space $X$ with topology $\tau$, the \emph{free topological group $F(X)$} (the \emph{free Abelian 
topological group $A(X)$}, the \emph{free Boolean topological group $B(X)$}) is the topological (topological 
Abelian, topological Boolean) group containing $X$ as a subspace, generated by $X$, and defined by the 
\emph{universal property} that any continuous function from $X$ to a topological group $G$ 
(topological Abelian group $G$, topological Boolean group $G$) extends to a continuous homomorphism 
$F(X)\to G$ ($A(X)\to G$, $B(X)\to G$). In other words, this is the abstract free (free 
Abelian, free Boolean) group of the set $X$ endowed with the finest group topology inducing the topology 
$\tau$ on $X$. Basic information about the groups $F(X)$ and $A(X)$ can be found in \cite{VINITI} and~\cite{AT}; 
the groups $B(X)$ were studied in~\cite{Axioms}. 

If $X$ is zero-dimensional, then the \emph{free linear topological group 
$F^{\mathrm{lin}}(X)$}  (as well as the \emph{free Abelian linear topological group $A^{\mathrm{lin}}(X)$} and 
the \emph{free Boolean linear topological group $B^{\mathrm{lin}}(X)$}) is also defined 
\cite[Theorem~2]{Axioms}. Its definition is similar to that of the free topological group of $X$, the only 
difference being that its topology is required to be linear and continuous functions from $X$ to $G$ must 
extend to continuous homomorphisms only for $G$ with linear topology. For more details concerning free linear 
topological groups, see~\cite{Axioms}. 

\begin{remark}
\label{r2}
Thanks to the universal property, $A(X)$ is a continuous homomorphic image of $F(X)$, $B(X)$ is a continuous 
homomorphic image of $A(X)$, and $F^{\mathrm{lin}}(X)$, $A^{\mathrm{lin}}(X)$, and 
$B^{\mathrm{lin}}(X)$ are the images of $F(X)$, $A(X)$, and $B(X)$, respectively, under the continuous identity 
isomorphisms. 
\end{remark}

The topology of any topological group $G$ with identity element $1$ 
is induced by the natural \emph{two-sided group uniformity} $\mathscr V_G$
with base 
$$ 
\{\{(g,h)\in G\times G: h\in gV\cap Vg\}: \text{$V$ is a neighborhood of $1$}\}
$$
(see \cite[Example~8.1.17]{Eng}). A topological group $G$ is said to be 
\emph{Raikov complete} if $\mathscr V_G$ is complete. It is well known that $G$ is Raikov complete if and only if 
it is closed in any topological group containing $G$ as a topological subgroup \cite{Raikov}. 
More details on Raikov complete groups can be found in \cite[Section~3.6]{AT} (see also~\cite{Raikov}). 

A topological group $G$ is said to be \emph{$\mathbb R$-factorizable} if, for every continuous function 
$f\colon G\to \mathbb R$, there exists a continuous homomorphism $h\colon G \to H$ to a second-countable 
topological group $H$ and a continuous function $g\colon H\to \mathbb R$ such that $f = g \circ h$. This 
very useful notion was introduced by Tkachenko \cite{Tk1}, who showed, among other things, that any Lindel\"of 
group is $\mathbb R$-factorizable~\cite[Assertion~1.1]{Tk1} and that any $\mathbb R$-factorizable group $G$ is 
\emph{$\omega$-narrow}, that is, for every neighborhood $U$ of the identity element in $G$, there 
exists a countable set $A\subset G$ for which $A\cdot U=G$ (see \cite[Proposition~8.1.3]{AT}).

In what follows, we repeatedly use the following known theorems. 

\begin{theoremA}[{see, e.g., \cite[Corollary~7.1.18]{AT}}]
\label{tA}
The free topological group $F(X)$ of a space $X$ is Lindel\"of if and only if $X^n$ is 
Lindel\"of for each $n \in \mathbb N$. 

Therefore, if $X^n$ is Lindel\"of for each $n \in \mathbb N$, then $A(X)$, $B(X)$, $F^{\mathrm{lin}}(X)$, 
$A^{\mathrm{lin}}(X)$, and $B^{\mathrm{lin}}(X)$ are Lindel\"of. 
\end{theoremA}

\begin{theoremA}[{\cite[Theorem~3.1]{Sh}; see also \cite[Theorem~8.8.4]{AT}}]
\label{tB}
Suppose that $G$ is a zero-dimensional $\mathbb R$-factorizable group, $H$ is a second-countable topological 
group, and $f\colon G \to H$ is a continuous homomorphism. Then there exists a zero-di\-men\-sional 
second-countable topological group $G'$ and continuous epimorphisms $g\colon G \to G'$ and $h\colon G' \to H$ 
such that $f = g \circ h$. 
\end{theoremA}

\begin{theoremA}[{\cite{Morita1}, \cite{Nagami}; see also \cite[Theorem~11.22]{Ch}}]
\label{tC}
If $Y$ is a $z$-embedded subspace of a space $X$, then $\dim_0 Y \le \dim_0 X$.
\end{theoremA}

The definitions and facts used in this paper without reference can be found in \cite{Eng} or~\cite{AT}.

\section{There is No Product Theorem\\ for the Covering Dimension of Topological Groups}

\begin{theorem}
\label{t1}
There exist Boolean (and hence Abelian) topological groups $G_1$ and $G_2$ with the following properties:
\begin{enumerate}
\item[\textup{(1)}]
$G_1^n$ is Lindel\"of and submetrizable for every $n\in \mathbb N$\textup;
\item[\textup{(2)}]
$G_2$ is second-countable;
\item[\textup{(3)}]
the topologies of $G_1$ and $G_2$ are linear;
\item[\textup{(4)}]
$\dim_0 G_1 = \dim G_1 = 0$ and $\dim_0 G_2 = \dim G_2 = 0$\textup;
\item[\textup{(5)}]
$\dim (G_1\times G_2)>0$ and $\dim_0 (G_1\times G_2)>0$. 
\end{enumerate}
\end{theorem}

To prove the theorem, we need two lemmas.

\begin{lemma}
\label{l1}
If a space $X$ is a retract of a topological group $G$, then it is a retract of the free topological 
group $F(X)$. If $G$ is Abelian \textup(Boolean\textup), then $X$ is a 
retract of the free Abelian topological group $A(X)$ \textup(of the free Boolean topological 
group~$B(X)$\textup). If the topology of $G$ is linear, then $X$ is a retract of the free linear topological 
group $F^{\mathrm{lin}}(X)$\textup; if, in addition, $G$ is Abelian \textup(Boolean\textup), then $X$ is a retract of 
$A^{\mathrm{lin}}(X)$  \textup(of $B^{\mathrm{lin}}(X)$\textup). 
\end{lemma}

\begin{proof}
Let $r\colon G\to X$ be a retraction. Then the restriction $r|_{\langle X\rangle}$ of $r$ to the subgroup 
$\langle X\rangle $ of $G$ generated by $X$ is a retraction as well, because $X\subset \langle X\rangle$. By the 
definition of $F(X)$ the identity map $\mathrm{id}_X\colon X\to X\subset \langle X\rangle$ extends to a 
continuous homomorphism $h\colon F(X)\to \langle X\rangle$. Clearly, $r|_{\langle X\rangle}\circ h$ is a 
retraction. 

In the cases where $G$ is Abelian or Boolean and where the topology of $G$ is linear, the argument is similar.
\end{proof}

\begin{lemma}
\label{l2}
Every second-countable space $X$ which is a retract of a topological group $G$ is a retract of 
a topological group $H$ with the following properties: 
\begin{enumerate}
\item[\textup{(1)}]
$H$ is second-countable;
\item[\textup{(2)}]
if $G$ is Abelian or Boolean, then so is~$H$;
\item[\textup{(3)}] 
if $X$ is zero-dimensional, then so is~$H$;
\item[\textup{(4)}] 
if $G$ is Abelian\footnotemark\ and its topology is linear, then so is the topology of~$H$.
\end{enumerate}
\end{lemma}

\footnotetext{This assumption is made for simplicity, it can be dropped.}

\begin{proof}
Suppose that a second-countable space $X$ is a retract of a topological group $G$. By Lemma~\ref{l1}\, $X$ is a 
retract of the free topological group $F(X)$; let $r$ be a retraction $F(X)\to X$. According to Theorem~\ref{tA}, 
$F(X)$ is Lindel\"of and hence $\mathbb R$-factorizable. By Assertion~1.1 of \cite{Tk1} there exists a 
second-countable group $H$, a continuous epimorphism $h\colon F(X)\to H$, and a continuous map $f\colon H\to X$ 
for which $r= f\circ h$. 

Let $\mathrm{id}_X$ denote the identity embedding of $X$ into $F(X)$, and let $f' = h\circ 
\mathrm{id}_X\colon X\to H$. Clearly, $f'$ is continuous and $f(f'(x))= f(h(x))= r(x)=x$ for 
$x\in X$. According to \cite[Theorem on p.~1085]{Borsuk}, $X$ is homeomorphic to a retract of~$H$. 

If $G$ is Abelian or Boolean, then we render $H$ Abelian or Boolean by replacing  $F(X)$ 
with $A(X)$ or $B(X)$, respectively, in the above argument. The groups $A(X)$ and $B(X)$ are Lindel\"of by 
Theorem~\ref{tA}.

If $X$ is zero-dimensional (and hence strongly zero-dimensional, being second-countable), 
then so are $F(X)$ \cite[Proposition~1]{Arh68} (see also \cite[Theorem~7.6.16]{AT}), $A(X)$ \cite{Tk-0}, 
and $B(X)$ \cite[Theorem~8]{Axioms}. Thus, in this case, Theorem~\ref{tB} applies, according to which 
the group $H$ can be made zero-dimensional. 

Suppose that $G$ is Abelian and its topology is linear. Then $X$ is a retract of the free Abelian linear 
topological group $A^{\mathrm{lin}}(X)$ (by Lemma~\ref{l1}). Let $r\colon A^{\mathrm{lin}}(X)\to X$ be a 
retraction. By Theorem~\ref{tA}\, $A^{\mathrm{lin}}(X)$ is Lindel\"of. As 
above, applying Assertion~1.1 of \cite{Tk1}, we find a second-countable group $\widetilde H$, a continuous 
epimorphism $\tilde h\colon A^{\mathrm{lin}}(X)\to \widetilde H$, and a continuous map $\tilde f\colon \widetilde 
H\to X$ for which $r= \tilde f\circ \tilde h$. Note that $\widetilde H$ is Abelian. Let $\{U_n: n\in \omega\}$ be 
a base of neighborhoods of zero in~$\widetilde H$. For each $n$, $\tilde h^{-1}(U_n)$ contains an open subgroup 
$A_n$ of $A^{\mathrm{lin}}(X)$; we set $H_n=\tilde h(A_n)$. The subgroups $H_n$ are normal and hence form a 
subbase of neighborhoods of zero for some group topology on $\widetilde H$; let $H$ be $\widetilde H$ with this 
new topology. The new topology is finer than the old one; therefore, the map $f\colon H\to X$ coinciding with 
$\tilde f$ as a map of sets is continuous. The epimorphism $h\colon A^{\mathrm{lin}}(X)\to H$ coinciding with 
$\tilde h$ as a map of abstract groups is continuous as well, because the preimage of any basic neighborhood 
of zero in $H$ contains an open neighborhood of zero in $A^{\mathrm{lin}}(X)$. The group $H$ is second-countable, 
because it is metrizable (being first-countable) and Lindel\"of (being a continuous image of the Lindel\"of group 
$A^{\mathrm{lin}}(X)$), and $X$ is a retract of $H$ from the same considerations as in the second paragraph of 
this proof. 

If $G$ is Boolean, then $H$ can be rendered Boolean by considering $B^{\mathrm{lin}}(X)$  
instead of $A^{\mathrm{lin}}(X)$.
\end{proof}

\begin{proof}[Proof of Theorem~\ref{t1}]
Our construction of the topological groups $G_1$ and $G_2$ is based on Charalambous' 
modification of Przymusi\'nski's construction in~\cite{Prz} of a strongly zero-dimensional Lindel\"of space 
whose square is normal but not strongly zero-dimensional. Namely, in the proof of Theorem~27.5 in \cite{Ch} 
subsets $S$, $S_1$ and $S_2$ of the Cantor set $C$ and topologies $\tau_1$ and $\tau_2$ on $C$ with certain 
properties were defined. We put $C_1=(C,\tau_1)$ and denote the usual Euclidean topology of $C$ by 
$\tau$. The set $S_2$ is assumed to be endowed with the topology induced by $\tau_2$, which coincides 
with that induced by $\tau$ (that is, by the usual topology of~$C$). 

We need the following properties of $\tau_1$,  $C_1$ and $S_2$: 
\begin{enumerate} 
\item[\textup{(i)}] 
$\tau_1$ is finer than $\tau$; 
\item[\textup{(ii)}] 
$\tau_1$ has a base consisting of sets closed in $\tau$; 
\item[\textup{(iii)}] 
$C_1$ is first-countable and Lindel\"of; 
\item[\textup{(iv)}] 
$S_2$ is second-countable;
\item[\textup{(v)}] 
$\dim C_1=\dim_0 C_1=0$; 
\item[\textup{(vi)}] 
$\dim S_2=\dim_0 S_2=0$.  
\end{enumerate}
 In \cite[Example~27.8]{Ch} it was shown that 
\begin{enumerate} 
\item[\textup{(vii)}] 
$\dim_0(C_1\times S_2)>0$ (and hence $\dim (C_1\times S_2)\ge 0$). 
\end{enumerate} 
In \cite{arx} the construction was refined so as to satisfy the additional condition 
\begin{enumerate} 
\item[\textup{(viii)}] 
$C_1^n$ is Lindel\"of for each $n\in \mathbb N$. 
\end{enumerate}

Recall that a topological space is said to be \emph{non-Archimedean} if it has a base such that, given  
any two of its elements, either they are disjoint or one of them contains the other (see~\cite{Nyikos}). Note 
that the Cantor set $C$ (with the usual topology), as well as its subspace $S_2$, is non-Archimedean. 
According to \cite[Theorem~3 (version~2)]{3}, any space $X$ admitting a coarser non-Archimedean topology 
$\sigma$ and having a base consisting of $\sigma$-closed sets is a retract of a Boolean topological group 
with linear topology. By Lemma~\ref{l1} any such $X$ is a retract of $B^{\mathrm{lin}}(X)$ (in fact, the group 
constructed in \cite{3} \emph{is} $B^{\mathrm{lin}}(X)$). Thus,  $C_1$ and $S_2$ are retracts of the 
zero-dimensional Boolean groups $B^{\mathrm{lin}}(C_1)$ and $B^{\mathrm{lin}}(S_2)$, respectively. 

For each $n\in \mathbb N$, we denote the topological sum of $n$ copies of $C_1$ by $\bigoplus_n C_1$. According 
to Proposition~7 of \cite{Axioms}, $(B(C_1))^n$ is topologically isomorphic to the group 
$B(\bigoplus_n C_1)$, which is Lindel\"of by Theorem~\ref{tA}. Since $B^{\mathrm{lin}}(C_1)$ is a continuous 
image of $B(C_1)$, it follows that all finite powers $(B^{\mathrm{lin}}(C_1))^n$ are Lindel\"of. 

Let us show that $B^{\mathrm{lin}}(C_1)$ is submetrizable. Since $C_1=(C,\tau_1)$ and the topology $\tau_1$ 
is finer than the Euclidean topology $\tau$ of the Cantor set $C$, it follows that the 
identity isomorphism $B^{\mathrm{lin}}(C_1)\to B^{\mathrm{lin}}(C)$ extending the identity map $C_1\to C$ is 
continuous. On the other hand, $B^{\mathrm{lin}}(C)$ is a continuous image of $F(C)$ and $F(C)$ has a countable 
network \cite[Theorem~5.2.13]{AT}; therefore, $B^{\mathrm{lin}}(C)$ has a countable network as well and hence 
admits a coarser metrizable group topology~\cite[Corollary~7.1.17]{AT}, which immediately implies the 
submetrizability of $B^{\mathrm{lin}}(C_1)$.

We set $G_1=B^{\mathrm{lin}}(C_1)$ and let $r_1$ be a retraction $G_1\to C_1$. 

By Lemma~\ref{l2}\, $S_2$ is a retract of a Boolean second-countable topological group with linear 
topology. We denote this group by $G_2$ and let $r_2$ be a retraction $G_2\to S_2$.  

Since the dimension $\dim$ of a Lindel\"of space does not exceed its small inductive dimension $\ind$ (see, 
e.g., \cite[Proposition~5.3]{Ch}) and the dimensions $\dim_0$ and $\dim$ coincide for normal spaces (see, e.g., 
\cite[Proposition~11.2]{Ch}), it follows that $\dim G_1=\dim_0 G_1=0$ and $\dim G_2=\dim_0 G_2=0$. However, 
$\dim_0 (G_1\times G_2) > 0$. Indeed, clearly, $r_1\times r_2\colon G_1\times G_2 \to C_1\times S_2$ is a 
retraction. Thus, $C_1\times S_2$ is a retract and hence a $z$-embedded subspace of $G_1\times G_2$. According to 
Theorem~\ref{tC}, in view of (vii) and Remark~\ref{r1} we have $\dim_0 (G_1\times G_2)>0$ and $\dim (G_1\times 
G_2) > 0$. 
\end{proof}

\section{Covering Dimension is Not Preserved by Free Topological Groups}

\begin{theorem}
\label{th3}
There exists a space $X$ with the following properties:
\begin{enumerate}
\item[(1)] 
$\dim X=\dim_0 X=0$\textup;
\item[(2)]
$\dim F(X)>0$ and $\dim_0 F(X)>0$\textup;
\item[(3)]
$\dim A(X) >0$ and $\dim_0 A(X)>0$\textup;
\item[(4)]
$\dim B(X) >0$ and $\dim_0 B(X)>0$.
\end{enumerate}
\end{theorem}

\begin{proof}
We keep the notation of the preceding section. Let us show that 
$X=C_1\oplus S_2$ has the desired properties. 

Property (1) obviously follows from properties (v) and (vi) of the spaces $C_1$ and $S_2$ (see the proof of 
Theorem~\ref{t1}). Properties (3) and (4) follow from property (vii) of $C_1\times S_2$ and 
the fact that, according to \cite[Proposition~4]{Tkachuk} and \cite[Proposition~7]{Axioms}, the groups 
$A(C_1)\times A(S_2)$ and $B(C_1)\times B(S_2)$ are topologically isomorphic to $A(C_1\oplus S_2)$ and 
$B(C_1\oplus S_2)$, respectively. Indeed, we know from the proof of Theorem~\ref{t1} that $C_1$ and $S_2$ are 
retracts of Boolean topological groups. By Lemma~\ref{l1} they are also retracts of $A(C_1)$ and 
$B(C_1)$ and of $A(S_2)$ and $B(S_2)$, respectively. Therefore, $C_1\times S_2$ is a retract of 
$A(C_1)\times A(S_2)$ and of $B(C_1)\times B(S_2)$. Thus, the groups $A(C_1\oplus S_2)\cong A(C_1)\times A(S_2)$ 
and $B(C_1\oplus S_2)\cong B(C_1)\times B(S_2)$ cannot be strongly zero-dimensional, because they contain the 
space $C_1\times S_2$ with $\dim_0(C_1\times S_2)>0$ as a $z$-embedded subspace. By Remark~\ref{r1} their 
dimension $\dim$ cannot be zero either. 

Let us prove (2). The natural multiplication map $i_2\colon X\times X\to F(X)$ defined by $(x,y)\mapsto xy$ is a 
topological embedding (see, e.g., \cite[Theorem~7.1.13]{AT}). Therefore, $C_1\times 
S_2$ is topologically embedded in $F(X)$ as the subspace $Y=i_2(C_1\times S_2)$ consisting of two-letter words of 
the form $xy$, where $x\in C_1$ and $y\in S_2$. Let us show that $Y$ is a retract of~$F(X)$.  

The group $A(X)=A(C_1\oplus S_2)$ is the topological quotient of $F(X)=F(C_1\oplus S_2)$ by the commutator 
subgroup (see, e.g., \cite[Theorem~7.1.11]{AT}). Let $h\colon F(C_1\oplus S_2)\to A(C_1\oplus S_2)$ be the 
canonical quotient homomorphism. Note that
$$ 
h(Y)=\{x+y\in A(C_1\oplus S_2): x\in C_1,\ y\in S_2\}. 
$$ 

The isomorphism $i\colon A(C_1\oplus S_2)\to A(C_1)\times A(S_2)$ constructed in \cite{Tkachuk}
takes each point $x\in C_1\oplus S_2$ to $(x,0_2)$ if $x\in C_1$ and to $(0_1, x)$ if $x\in S_2$ (by $0_1$ and 
$0_2$ we denote the zero elements of $A(C_1)$ and $A(S_2)$, respectively), so that 
$$ 
i(x+y)=(x,y)\in C_1\times S_2\subset A(C_1)\times A(S_2) 
$$ 
for any $x\in C_1$ and $y\in S_2$. Obviously, 
$$ 
i(h(Y))=C_1\times S_2\subset A(C_1)\times A(S_2). 
$$ 

As mentioned above, $C_1\times S_2$ is a retract of $A(C_1)\times A(S_2)$. Let $r\colon A(C_1)\times A(S_2) \to 
C_1\times S_2$ be a retraction. The composition $r\circ i \circ h\colon F(C_1\oplus S_2)\to C_1\times S_2$ is 
surjective and continuous,  and it takes every $xy\in Y$ to $(x,y)\in C_1\times S_2$. To obtain the desired 
retraction $F(X)\to Y$, it remains to add the homeomorphism $i_2|_{C_1\times S_2}\colon C_1\times S_2\to Y$  to 
this composition.  

Thus, $Y$ is a retract and hence a $z$-embedded subspace of $F(X)$. Since $Y$ is homeomorphic to $C_1\times 
S_2$, we have $\dim_0 Y>0$ (see property (vii) of $C_1\times S_1$). Therefore, $\dim_0 
F(X)>0$ by Theorem~\ref{tC} and $\dim F(X)>0$ by Remark~\ref{r1}. 
\end{proof}

\section{There is No Subgroup Theorem\\ for the Covering Dimension of Topological Groups}

\begin{theorem}
\label{t2}
There exists a strongly zero-dimensional Boolean (and hence Abelian) group $G$ topology which 
contains a closed subgroup $H$ with $\dim_0 H>0$. 
\end{theorem}

The proof of this theorem uses the following lemma.

\begin{lemma}
\label{l3}
Any zero-dimensional $\mathbb R$-factorizable group $G$ embeds in a product $P$ of zero-dimensional 
second-countable groups as a subgroup. Moreover, if $G$ is Abelian or Boolean, then so is~$P$, and if in 
addition $G$ has linear topology, then so does $P$. 
\end{lemma}

\begin{proof}
Let $G$ be a zero-dimensional $\mathbb R$-factorizable group, and let $f_\alpha\colon G\to \mathbb R$, $\alpha 
\in A$, be all continuous functions on $G$ (here $A$ is some index set). 
It follows from the $\mathbb R$-factorizability of $G$ and Theorem~\ref{tB} that, for each $\alpha\in A$, there 
exists a zero-dimensional second-countable group $H_\alpha$, a continuous epimorphism $h_\alpha\colon G\to 
H_\alpha$, and a continuous function $g_\alpha\colon H_\alpha \to \mathbb R$ such that $f_\alpha= g_\alpha\circ 
h_\alpha$. Note that if $G$ is Abelian or Boolean, then so are all $H_\alpha$. Since $G$ is Tychonoff, it follows 
that the family $\{f_\alpha:\alpha \in A\}$ separates points and closed sets and hence so does $\{h_\alpha:\alpha 
\in A\}$. Therefore, the diagonal 
$$ 
\mathop{\raise-3pt\hbox{\LARGE $\Delta$}}\limits_{\alpha \in A} h_\alpha\colon G\to \prod_{\alpha\in A} H_\alpha
$$ 
is a homeomorphic embedding. Clearly, this is a homomorphism. We set $P=\prod_{\alpha\in A} H_\alpha$.

In the case where $G$ is Abelian and has linear topology, we can render the topologies of all $H_\alpha$ linear 
in the same manner as in the proof of Lemma~\ref{l2}: for each $\alpha\in A$, we fix a base $\{U_n: n\in 
\omega\}$ of neighborhoods of zero in~$H_\alpha$, choose a subgroup of $H_\alpha$ with open preimage under 
$h_\alpha$ in each $U_n$, and define the new group topology on $H_\alpha$ for which the chosen subgroups form a 
subbase of neighborhoods of zero. The maps $h_\alpha$ and $f_\alpha$ remain continuous with respect to the new 
topology, and the group $H$ with this topology is first-countable and $\omega$-narrow, because $G$ is 
$\omega$-narrow, being $\mathbb R$-factorizable, and continuous homomorphisms preserve $\omega$-narrowness 
\cite[Proposition~3.4.2]{AT}. Therefore, it is second-countable \cite[Proposition~3.4.5]{AT}. 

Clearly, the topology of any product of groups with linear topology is linear. Therefore, the product $P$ of the 
second-countable groups $H_\alpha$ with the new linear topologies is linear. 
\end{proof}

\begin{proof}[Proof of Theorem~\ref{t2}]
We use the same spaces $C_1$ and $S_2$ as in the proof of Theorem~\ref{t1}. By Lemma~\ref{l1}\, $C_1$ and $S_2$ 
are retracts of the free Boolean groups $B(C_1)$ and $B(S_2)$, respectively. These groups are Lindel\"of 
by Theorem~\ref{tA}. According to \cite[Theorem~8]{Axioms}, they are zero-dimensional, and according to 
\cite{Tk1}, they are $\mathbb R$-factorizable. By Lemma~\ref{l3}\, $B(C_1)$ and $B(S_2)$ are embedded in products 
$P_1$ and $P_2$ of zero-dimensional second-countable Boolean groups as subgroups. By Theorem~3 of \cite{Morita} 
any product of zero-dimensional second-countable spaces is strongly zero-dimensional. Therefore, the group 
$G=P_1\times P_2$ is strongly zero-dimensional, and it contains $H=B(C_1)\times B(S_2)$ as a subgroup. By 
Theorem~\ref{tC}\, $\dim_0 H>0$, because $\dim_0(C_1\times S_2)>0$ and $C_1\times S_2$ is a retract and hence a 
$z$-embedded subspace of~$H$. 

The subgroup $H$ is closed in $G$, because it is Raikov complete. Indeed, 
$C_1$ and $S_2$ are paracompact, being Lindel\"of, and therefore Dieudonn\'e complete~\cite{Dickinson}. It 
follows from Theorem~2.1 of \cite{arx2} that the free Boolean topological group of any Dieudonn\'e complete space 
is Raikov complete. It remains to recall that Raikov completeness is preserved by products (see, e.g., 
\cite[Theorem~3.6.22]{AT}). 
\end{proof}

\begin{remark}
It is easy to construct a similar example for the \v Cech covering dimension 
$\dim$ (but the subgroup $H$ cannot be made closed in this case, because, for the dimension $\dim$, the closed 
subset theorem holds \cite[Proposition~2.11]{Ch}). Indeed, consider the Sorgenfrey plane $S\times S$. It is 
zero-dimensional and has weight $2^\omega$. Hence it embeds in the Cantor cube $K=\{0,1\}^{2^\omega}$ 
\cite[Theorem~6.2.16]{Eng}. Since $K$ is strongly zero-dimensional, it follows that the free Abelian topological 
group $A(K)$ is zero-dimensional \cite{Tk-0}, and since $A(K)$ is Lindel\"of (by Theorem~\ref{tA}), it follows 
 that $\dim A(K)=0$ \cite[Proposition~5.3]{Ch}. Let $H$ denote the subgroup of $A(K)$ generated by 
$S\times S$. Clearly, $H\cap K=S\times S$, and hence $S\times S$ is closed in $H$. Therefore, $\dim H=\infty$, 
because $\dim (S\times S)=\infty$~\cite{Dijkstra}. 
\end{remark}

\begin{remark}
Applying the argument of the proof of Theorem~\ref{t2} to the free Boolean linear topological groups 
$B^{\mathrm{lin}}(C_1)$ and $B^{\mathrm{lin}}(S_2)$ instead of $B(C_1)$ and $B(S_2)$, we obtain an example of a 
strongly zero-dimensional Boolean group $G$ with linear topology which contains a subgroup $H$ with $\dim_0 
H>0$. However, it is unclear whether $H$ can be made closed, because the free Boolean linear group of a 
Dieudonn\'e complete (and even compact) space is not necessarily complete. For example, a base of neighborhoods 
of zero in the group $B^{\mathrm{lin}}(\xi)=[\xi]^{<\omega}$ for the usual convergent sequence $\xi=\mathbb N\cup 
\{\infty\}$ is formed by the subgroups 
$$ 
H_n=\{F\subset \xi\setminus \{1,\dots, n\}: \text{$|F|$ is even}\}
$$
(see the description of the topology of $B^{\mathrm{lin}}(X)$ in \cite[p.~497]{Axioms}). It is easy to see that 
$B^{\mathrm{lin}}(\xi)$ is topologically isomorphic to the $\sigma$-product of countably many copies of the 
discrete group $\mathbb Z_2=\{0,1\}$: the isomorphism takes every element $F\in B^{\mathrm{lin}}(\xi)$ (which is 
a finite subset of $\xi$) to the point $(x_n)_{n\in \omega}\in \{0,1\}$ in which $x_0=1$ if and only if 
$\infty \in F$ and $x_n=1$ for $n>0$ if and only if $n\in F$. This $\sigma$-product is a dense and hence 
nonclosed subgroup of $\mathbb Z_2^\omega$. Therefore, it is not complete. 
\end{remark}

The author thanks Evgenii Reznichenko for discussions.


\begin{thebibliography}{00}

\bibitem{Arh68}
A. V. Arhangel'skii, 
``Mappings connected with topological groups,''
Dokl. Akad. Nauk SSSR \textbf{181}, 1303--1306 (1968).

\bibitem{Arh89}
A.~V.~Arkhangelskii, 
``Algebraic objects generated by topological structure,''
J. Math. Sci. (N.Y.) \textbf{45} (1), 956--990 (1989). 

\bibitem{A-vM}
A.~V.~Arhangel'skii and J.~van~Mill, 
``Some aspects of dimension theory for topological
groups,'' Indag. Math. \textbf{29}, 202--225 (2018).

\bibitem{AT}
A.~Arhangel'skii and M.~Tkachenko,
\textit{Topological Groups
and Related Structures}
(Atlantis Press/World Sci., Amsterdam--Paris, 2008).

\bibitem{Borsuk}
K. Borsuk, 
``On The Topology of Retracts,'' 
Ann. Math. \textbf{48} (4), 1082--1094 (1947). 

\bibitem{Ch}
M. G. Charalambous,
\textit{Dimension Theory:
A Selection of Theorems
and Counterexamples} 
(Springer International, Cham, 2019).

\bibitem{Dickinson}
A. Dickinson, 
``Compactness conditions and uniform structures,'' 
Amer. J. Math. \textbf{75} (2), 224--228 (1953).

\bibitem{Dijkstra}
J. J. Dijkstra, 
``A space with maximal discrepancy between its Kat\v etov and covering dimension,''
Topol. Appl. \textbf{12} (1), 45--48 (1981).

\bibitem{Eng}
R. Engelking, 
\textit{General Topology}, 2nd ed. (Heldermann-Verlag, Berlin, 1989).

\bibitem{3}
P.~M.~Gartside, E.~A.~Reznichenko, and O.~V.~Sipacheva, 
``Mal'tsev and retral spaces,''
Topol. Appl. \textbf{80}, 115--129  (1997).

\bibitem{Morita}
 K. Morita, 
``On the dimension of the product of Tychonoff spaces,'' General Topology Appl. \textbf{3},
123--133 (1973).

\bibitem{Morita1}
 K. Morita, 
``Dimension of General Topological Spaces,'' in: \textit{Surveys in General Topology} (Academic Press, New York, 
1980), pp.~297--336. 

\bibitem{Nagami}
 K. Nagami, 
``Dimension of non-normal spaces,'' 
Fund. Math. \textbf{109}, 113--121 (1980).

\bibitem{Nyikos}
P. J. Nyikos, 
``Some surprising base properties in topology,''
in \textit{Studies in Topology}, ed.\ by  N.~M.~Stavrakas and K.~R.~Allen (Academic Press, 
New York--San~Francisco--London, 1975),
pp.~427--450. 

\bibitem{Prz}
T.~C.~Przymusi\'nski, 
``On the dimension of product spaces and an example of M.~Wage,'' Proc. Amer.
Math. Soc. \textbf{76}, 315--321 (1979).

\bibitem{Raikov}
D. Raikov, 
``On the completion of topological groups,'' 
Izv. Akad. Nauk SSSR, Ser. Mat. \textbf{10}, 513--528 (1946).

\bibitem{Sh}
D. B. Shakhmatov, 
``Imbeddings into topological groups preserving dimensions,''
Topol. Appl. \textbf{36}, 181--204 (1990).

\bibitem{88}
D. Shakhmatov, 
``A survey of current researches and open problems in the dimension theory of topological
groups,'' Questions Answers Gen. Topology \textbf{8} (1), 101--128  (1990).

\bibitem{VINITI}
O. V. Sipacheva, 
``The topology of free topological groups,'' 
J. Math. Sci. \textbf{131}, 5765--5838 (2005).

\bibitem{Axioms}
O. Sipacheva, 
``Free Boolean topological  groups,''
Axioms \textbf{4}, 492--517  (2015).

\bibitem{arx}
O. Sipacheva,
``A non-$\mathbb R$-factorizable product of $\mathbb R$-factorizable groups,''
\href{https://arxiv.org/abs/2303.08878v2}{arXiv:2303.08878v2 [math.GN]} (2025).

\bibitem{arx2}
O. V. Sipacheva and M. G. Tkachenko,
``The completeness of free Boolean topological groups,''
\href{https://arxiv.org/abs/2507.12289}{arXiv:2507.12289 [math.GN]} (2025).

\bibitem{Tk-0}
M.~G.~Tkachenko, 
``Zero-dimensionality of free topological groups,'' 
C. R. Acad. Bulgare Sci. \textbf{38} (2), 173--174 (1985).

\bibitem{Tk1}
M.~G.~Tkachenko, 
``Factorization theorems for topological groups and their applications,'' 
Topol. Appl. \textbf{38}, 21--37 (1991). 

\bibitem{Tk}
Michael Tka\v cenko,
`` Subgroups, quotient groups and
products of $\bR$-factorizable groups,''
Topol. Proc. \textbf{16}, 201--231 (1991).

\bibitem{Tk4}
M. G. Tkachenko,
``Topological Features of Topological Groups,'' 
in \textit{Handbook of 
the History of General Topology}, Vol.~3, ed.\ by  C. E. Aull and R.  Lowen (Springer, Dordrecht, 2001), 
pp.~1027--1144. 

\bibitem{Tkachuk}
V.~V.~Tkachuk, 
``Duality with respect to the functor $C_p$ and cardinal invariants of the type of the Souslin number,''
Math. Notes \textbf{37}, 247--252 (1985). 

\end{thebibliography}
\end{document}